\documentstyle [epic,eepic,amssymb,amscd,11pt]{amsart}

\newtheorem{theorem}{Theorem}[section]
\newtheorem{Theorem}[theorem]{Theorem}
\newtheorem{definition}[theorem]{Definition}

\newtheorem{lemma}[theorem]{Lemma}
\newtheorem{Lemma}[theorem]{Lemma}
\newtheorem{corollary}[theorem]{Corollary}

\newtheorem{proposition}[theorem]{Proposition}
\newtheorem{example}[theorem]{Example}

\newtheorem{remark}[theorem]{Remark}

\newcommand{\pr} {\smallskip\noindent{\bf Proof\,\,}}
\newenvironment{Proof}	{\pr}{\hspace*{\fill}\qed\\}
\newenvironment{proof}	{\pr}{\hspace*{\fill}\qed\\}

\newcommand\R{{\Bbb R}}

\newcommand\C{{\Bbb C}}

\newcommand\Z{{\Bbb Z}}

\newcommand\cm{{\cal M}}
\newcommand\cC{{\cal C}}
\newcommand\cp{{\cal P}}

\newcommand\fg {{\frak g}}
\newcommand\fh {{\frak h}}
\newcommand\fk {{\frak k}}
\newcommand\ft {{\frak t}}

\newcommand\cS {{\cal S}}
\newcommand\U {{\cal U}}
\newcommand\tcU {{\tilde{\cal U}}}
\newcommand\tcV {{\tilde{\cal V}}}

\def	\tU	{{ \tilde{U}  }}
\def	\tG	{{ \tilde{G}  }}
\def	\tV	{{ \tilde{V}  }}
\def	\tT	{{ \tilde{T}  }}
\newcommand{\im}	{\operatorname{im}}
\newcommand{\Sl}	{\operatorname{SL}}

\begin{document}

\title{Symplectic toric orbifolds }
\author{Eugene Lerman  and Susan Tolman}\date{12/23/94}
\address{Department of Mathematics, MIT, Cambridge, MA
02139}
\email{eugene@@math.mit.edu, tolman@@math.mit.edu}
\thanks{Both authors were partially supported by NSF
postdoctoral fellowships}

\maketitle
\begin{center}
{ dg-ga/9412005}
\end{center}
\begin{abstract}
A  symplectic toric  orbifold is a compact connected orbifold $M$,
a symplectic form $\omega$ on $M$, and an effective
Hamiltonian action of a torus $T$ on $M$,
where the dimension of $T$ is half the dimension of $M$.

We prove that there is a one-to-one correspondence between symplectic
toric orbifolds and convex rational simple polytopes with positive
integers attached to each facet.
\end{abstract}

%\date{12/19/94}
%\input{intro}
\section{Introduction}
The main purpose of this paper is to demonstrate the following theorem.\\

\noindent
{\bf Theorem }
{\em Symplectic toric orbifolds
are classified by convex rational simple polytopes with a positive
integer attached to each facet.}\\

At the same time, however, we would like to build a
foundation for further work on symplectic orbifolds.  For this reason,
we will try to state lemmas in their most natural generality, rather
than restricting to our special case.

The above theorem  generalizes a theorem of Delzant \cite{Del}
to the case of orbifolds.
He proved that symplectic toric {\em  manifolds} are classified by
the image of their moment maps, that is,
by a certain class of rational polytopes.
It is easy to see that additional information is necessary in our case:
\begin{example}{\em
Given positive integers $r$ and $s$, let $l$ be their greatest common
divisor.
Let $K = \Z/(l\Z) \times S^1$ act on $\C^2$ by $(\xi,\lambda) \cdot (x,y)
= (\xi \lambda^n x, \lambda^m y)$, where $n = r/l$ and  $m = s/l$.
Let $(M,\omega)$ be the symplectic reduction of $\C^2$ with its standard
symplectic form at a positive number.
Then $T = (S^1)^2/K$ has an effective Hamiltonian action on $(M,\omega)$.
The image of the moment map is always a line interval, but
these spaces are not isomorphic.
Although  $M$ is (topologically) a two sphere,
it has two orbifold singularity, which look locally like $\C/(\Z/r\Z)$
and $\C/(\Z/s\Z)$.}
\end{example}

We begin by defining a few terms.

A {\em symplectic toric  orbifold} is a compact connected orbifold $M$,
a symplectic form $\omega$ on $M$, and an effective
Hamiltonian action of a torus $T$ on $M$,
where the dimension of $T$ is half the dimension of $M$.
This is not the definition used in algebraic geometry. We will see that every
symplectic toric orbifold can be given the structure of a projective toric
orbifold.
Two symplectic toric orbifolds are {\em isomorphic} if
they are equivariantly symplectomorphic.

Let $\ft$ be a vector space with a lattice $\ell$;
let $\ft^*$ denote the dual vector space.
%Note that the quotient $T =
%\ft/\ell $ is an $n$-torus with Lie algebra $\ft$.
A convex
polytope $\Delta \subset \ft^*$ is {\em rational} if the hyperplanes
supporting its facets are defined by the elements of the lattice
$\ell$, that is,
$$\Delta =
\cap_{i=1}^N \{\alpha \in \ft^* \mid \left< \alpha, y_i \right> \geq
\eta_i \}
$$
for some $y_i \in \ell$ and $\eta_i \in \R$. Recall that a {\em facet}
is a face of codimension one.  An $n$ dimensional polytope is {\em
simple} if exactly $n$ facets meet at every vertex.  For this paper,
we shall adopt the convenient but  non standard abbreviation that a {\em
weighted polytope} is a convex rational simple polytope plus a
positive integer attached to each facet.  Two
weighted polytopes are {\em isomorphic} if they differ by the
composition of a translation with an element of $\Sl(\ell) \equiv
\Sl(n,\Z)$ such that the corresponding facets have the same integers
attached to them.

Finally, let $x$ be a point in an orbifold $M$, and
let $(\tU,\Gamma, \phi)$ be a local uniformizing system for
neighborhood $U$ of $x$ (see \cite{Sa}),
then the {\em (orbifold) structure} group of $x$ is
the isotropy group of $\tilde{x} \in \tU$, where
$\phi(\tilde{x}) = x$.
This group is well defined.

We are now ready to give a precise statement of our main theorem:
\begin{theorem}
\label{classification}
For a symplectic toric orbifold $(M,\omega,T)$ the image of the moment
map is a rational simple polytope.
A positive integer $n$ is attached to every open facet of this polytope
as follows: for any $x$ in the preimage of the facet the
structure group of $x$ is $\Z/n\Z$.

Two symplectic toric orbifolds are isomorphic if and only
if their associated weighted polytopes are isomorphic.

Every weighted polytope occurs as the image of a
symplectic toric orbifold.
\end{theorem}

\begin{remark}{\em
We will see in section~\ref{surj} that all symplectic toric orbifolds
have K\"{a}hler structures.}
\end{remark}

\noindent {\sc Acknowledgments}\\
\vspace{-12pt}

At the conference {\em Applications of Symplectic Geometry} at the
Newton Institute, 10/31/94 - 11/11/94, we learned that R. de Souza and
E. Prato have been working independently on the same problem.

It is a pleasure to thank Chris Woodward for many useful discussions.
In particular, section~\ref{local to global} is joint work with Chris Woodward.
\\

\section{Group actions on orbifolds}

In this section, we recall facts about the action of groups
on orbifolds, following  Haefliger and Salem \cite{HS}.
Although some theorems about group actions on manifolds hold, there are a
few differences.
We begin with a few basic definitions.

An {\em orbifold} $M$ is a topological space $|M|$,
plus an {\em atlas} of {\em uniformizing charts}
$(\tU,\Gamma,\varphi)$, where
$\tU$ is open subset of $\R^n$, $\Gamma$ is a finite group which
acts linearly on $\tU$ and fixes a set of codimension at least two, and
$\varphi: \tU \to |M|$ induces a homeomorphism
from $\tU/\Gamma$ to $U \subset |M|$.
Just as for manifolds, these charts must cover $|M|$; they are also subject
to certain compatibility conditions; and there is a notion of when
two atlases of charts are equivalent.
For more details, see Satake \cite{Sa}.
Given $x$, we may choose a uniformizing chart $(\tU,\Gamma,\varphi)$
such that  $\varphi^{-1}(x)$ is a single point which is fixed by $\Gamma$.

Given $x \in M$, denote the {\em uniformizing tangent space at $x$} by
$\tT_xM$; it is $T_{\varphi^{-1}(x)}\tU$: the tangent space to
$\varphi^{-1}(x)$ in $\tU$.  Then $T_x M$, the fiber of the tangent
bundle of $M$ at $x$ is $\tT_x M/\Gamma$.  A vector field $\xi$ on $M$
is a $\Gamma$ invariant vector field $\tilde{\xi}$ on each
uniformizing chart $(\tU,\Gamma,\phi)$; of course, these must agree on
overlaps.  Similar definition apply to differential forms, etc.

Let $M$ and $N$ be orbifolds with atlases $\tcU$ and $\tcV$.
A {\em map of orbifolds} $f: M \to N$ is a map $f : \tcU \to \tcV$,
and an equivariant map $\tilde{f}_\tU: (\tU,\Gamma)
\to (\tV,\Upsilon)$ for each $(\tU,\Gamma) \in \tcU$,
where $(\tV,\Upsilon) = f(\tU,\Gamma)$.  These $\tilde{f}_\tU$
are subject to a compatibility condition
which insures, for instance, that $f$ induces a continuous map of
the underlying spaces.  Additionally, there is a notion of
when two such maps are equivalent.  Again, see \cite{Sa} for details.

We are now ready to recall the definition of the group action on an
orbifold.
\begin{definition}{\rm
Let $G$ be a Lie group.  A {\em smooth action} $a$ of $G$ on an
orbifold $M$ is a smooth orbifold map $a: G\times M \to M$ satisfying
the usual group laws, that is, for all $g_1, g_2 \in G$ and $x \in M$
$$
a(g_1, a(g_2, x )) = a(g_1 g_2, x ) \quad \hbox{and} \quad
a(1_G, x) = x,
$$
where $1_G$ is the identity element of $G$.  }
\end{definition}
Technically, by ``=", we mean ``are equivalent as maps of orbifolds."
The action $a$ induces a continuous action $|a|$ of $G$ on the
underlying topological space $|M|$,
$$
	|a|: G\times |M| \to |M|.
$$

In particular, this definition states that for every  $g_0 \in  G$
and $x_0 \in M$ there are neighborhoods $W$ of $g_0$ in $G$, $U$ of $x_0$
in $M$ and $U'$ of $a(g_0,x_0) $ in $M$, charts
$(\tilde{U}, \Gamma, \varphi) $ and $(\tilde{U}', \Gamma ', \varphi') $ and
a smooth map $\tilde{a}: W\times \tilde{U} \to \tilde{U}'$ such that
$\varphi' (\tilde{a}(g, \tilde{x})) =|a|(g, \varphi (x))$ for all $(g,
\tilde{x}) \in W\times \tilde{U}$.  Note that $\tilde{a}$ is not
unique, it is defined up to composition with elements of the orbifold
structure groups $\Gamma $ of $x_0$ and $\Gamma' $ of $g_0\cdot x_0$.

If $g_0 = 1_G$, the identity of $G$, then we may assume
$\tU \subset \tU'$, and we can choose $\tilde{a}$ such
that $\tilde{a}(1_G, x) = x$. Then $\tilde{a}$ induces a local action of
$G$ on $\tU$.  If in addition $x_0$ is {\em
fixed} by the action of $G$ and $G$ is compact, then the local action
$G$ on $\tU$ generates an action of $\tG$ on $\tV \subset \tU$,
where $\tilde{G}$ is a cover of the identity component of $G$.  Note that the
actions of $\tilde{G}$ and $\Gamma $ on $\tV$ {\em commute}.
%otherwise the group actions would not induce
%infinitesimal actions of their Lie algebras.

More generally one can show that for a fixed point $x$ with structure
group $\Gamma$ there exists
a uniformizing chart $(\tU,\Gamma,\phi)$ for a neighborhood $U$ of $x$,
an exact sequence of groups
$$
	1\to \Gamma \to \hat{G} \stackrel{\pi}{\to} G \to 1,
$$
%where the identity component  of $\hat{G}$ is $\tilde{G}$
and an action of $\hat{G}$ on $\tU$ such that the following diagram
commutes:
$$
\begin{CD}
\hat{G} \times \tU @>>> \tU\\
@VVV			@VVV\\
G \times U 	@>>>	U
\end{CD}
$$
The extension
$\hat{G}$ of $G$ depends on $x$ and, in particular, is not globally
defined.

For any $\xi \in \fg$, there is an associated vector field $\xi_M$
on $M$.  On each $\tU$, it is defined via the local action of
$G$ of $\tU$.
The map $\pi: \hat{G} \to G$ induces an isomorphism of Lie algebras
$\pi: \hat{\fg} \to  \fg$.
Locally, for any  $\hat{\xi} \in \hat{\fg}$,
$\hat{\xi}_\tU = \pi(\hat{\xi})_M |_\tU$.

\begin{remark}{\em
If $\hat\ft$ is a Cartan subalgebra of
$\hat{\fg}$, then $\ft = \pi(\hat{\ft})$ is also a Cartan subalgebra.
Define a lattice $\hat{\ell}$ by
$\hat{\ell} = \{ \hat{\xi} \in \hat{\ft} \mid \exp(\hat{\xi}) = 1\}$,
and define  $\ell \subset \ft$ analogously.
Then $\pi: \hat{\ell} \to \ell$ is not an isomorphism;
$\pi(\hat{\xi})$ is in $\ell$ exactly if $exp(\hat{\xi})$ is in $\Gamma$.
Let $\hat{\alpha}_i \in \hat{\ell}^*$ for $i \in I$
be the weights for the action of $\hat{G}$ on $\tilde{U}$.
Although $\alpha_i = \pi(\hat{\alpha}_i) \in \ft$ may not lie in the weight
lattice $\ell^*$,  $|\Gamma|  \alpha_i$ does lie in $\ell^*$, where
$|\Gamma|$ is the order of $\Gamma$.
We will call $\alpha_i$ the {\em orbi-weights} for the action of $G$ on $U$.
}\end{remark}

If $G$ is a compact  Lie group acting on an orbifold $M$,
we can define local slices to orbits in the differential category.
Since the orbit $G\cdot x$ is a {\em submanifold} of $M$,
$\tT_x (G\cdot x)$ is a $\Gamma$ invariant subspace of $\tT_x M$.
We define the {\em slice} at
$x$ for the action of $G$ on $M$ to be the orbifold $W/\Gamma$, where
$W= \tT_x M / \tT_x (G\cdot x)$.
We may also identify $W$
with the orthogonal complement to $\tT_x (G\cdot x)$ in $\tT_x M$ with
respect to some invariant metric.  It is not hard then to see that the
slice theorem takes the following form.

\begin{proposition}{\rm ({\bf Slice theorem})}
Suppose a compact Lie group $G$ acts on an orbifold $M$ and $G\cdot x$
is an orbit of $G$. Then a $G$ invariant neighborhood of the orbit is
equivariantly diffeomorphic to a neighborhood of the zero section in
the associated orbi-bundle $G\times _{G_x} W/\Gamma$, where $G_x$ is
the isotropy group of $x$ with respect to the action of $G$, $\Gamma$
is the orbifold structure group of $x$ and $W = \tT_x M/ T_x (G\cdot x)$
\end{proposition}

\begin{proof}
This is completely analogous to the slice theorem for actions on manifolds,
and follows immediately from the fact that metrics can be averaged over
compact Lie groups.
\end{proof}\vspace{-5mm}
\begin{remark}{\em
As in the smooth case, the compactness of the group $G$ is not necessary
for the existence of slices.  It is enough to require that the induced action
on the underling topological space is proper.}
\end{remark}

For connected $G$, it follows from the existence of slices
that the fixed point set is a suborbifold.
Therefore, the decomposition of M into infinitesimal orbit types is
a stratification into suborbifolds.
On the other hand, $M^G$ need not be a suborbifold in general.

\begin{example}{\em
Let $\Gamma = \Z/(2\Z)$ act on $\C^2$ by sending $(x,y)$ to $(-x,-y)$
Let $G = \Z/(2\Z)$ act on $\C^2/\Gamma$ by sending $[x,y]$ to $[x,-y]$.
Then $M^G = \{[x,0]\} \cup \{[0,y]\} = \C/\Gamma \cup \C/\Gamma$.}
\end{example}
Consequently the decomposition of an orbifold according
to the orbit type is not a stratification.
Fortunately the following lemma still holds.

\begin{lemma}\label{lem.princ.orbittype}
If $G$ is a compact Lie group acting on a connected orbifold $M$ then there
exists an open dense subset of $M$ consisting of points with the same orbit
type.
\end{lemma}
\begin{proof}
We first decompose the orbifold into the open dense set of smooth
points $M_{\text{smooth}}$ and the set of singular points. Since we
assume that all the singularities have codimension 2 or greater,
$M_{\text{smooth}}$ is connected. A smooth group action preserves this
decomposition. Since $G$ is compact, the action of $G$ on
$M_{\text{smooth}}$ has a principal orbit type (see for example
Theorem~4.27 in \cite{Kaw}). The set of points of this orbit type is
open and dense in $M_{\text{smooth}}$, hence open and dense in $M$.
\end{proof}

\begin{corollary}\label{cor.loc_free}
If a torus $T$ acts effectively on a connected orbifold $M$ then the action
of $T$ is free on a dense open subset of $M$.
\end{corollary}

\section{Symplectic local normal forms}
In this section, we write down normal forms for the
neighborhoods of orbits of compact Lie groups $G$ acting
symplectically on $(M,\omega)$;
that is, we classify such neighborhoods up to equivariant
symplectomorphisms.  We also point out consequences of this form,
both generally and in the special case of symplectic toric orbifolds.

A {\em symplectic orbifold} is an orbifold $M$ with a closed
nondegenerate $2$-form $\omega$.  A group $G$ acts {\em
symplectically} on $M$ if $G$ preserves $\omega$.  A moment map $\phi:
M \to \fg^*$ is a map such that for each $\tU$, the map $\tilde{\phi}$
is a moment map for the local group action of $G$ on $\tU$.  If there
is a moment map for $G$, we say that $M$ is a {\em Hamiltonian $G$
space}.

If  $G$ is a compact Lie group which acts symplectically
on $(M,\omega)$, we can define the {\em symplectic slice} at a point $x$.
The $2$-form $\omega$ induces a non-degenerate antisymmetric bilinear
form $\omega$ on $\tT_xM$.
Let $\tT (G\cdot x)^\omega $  be the symplectic perpendicular
to the tangent space  of $G \cdot x$ with respect to $\omega$.
The quotient
$$
V = \tT (G\cdot x)^{\omega }/( \tT (G\cdot x)\cap \tT (G\cdot x)^{\omega })
$$
is naturally a symplectic vector space.
The {\em symplectic slice } at $x$ is the symplectic orbifold $V/\Gamma,$
where $\Gamma$ is the structure group of $x$.

Notice that $\tT(G \cdot x)$ and $\tT(G \cdot x)^\omega$
are both $\hat{G}_x$ invariant, where $\pi: \hat{G}_x \to G_x$ is
an extension of $G_x$ by $\Gamma$.
Therefore there is a {\em symplectic linear action} of $G_x$ on $V/\Gamma$,
that is, a symplectic linear action of $\hat{G}_x$ on $V$.

\begin{remark} {\em
\label{rem_moment}
Let $\pi: \hat{\fg}_x \to \fg_x$ be the induced map of Lie algebras.
Let $\hat{\phi}_V: V \to \hat{\fg}_x^*$ be the moment map for  the
action of $\hat{G}_x$ on $V$.  Then
$\phi_{V/\Gamma}: V/\Gamma \to  \fg_x^*$, the moment
map  for the action of $G_x$ on $V/\Gamma$,
is given by the following diagram:
$$
\begin{CD}
V @>>> V/\Gamma \\
@V{\hat{\phi}_V}VV			@VV{\phi_{V/\Gamma}}V\\
\hat{\fg}_x^* 	@<\pi^*<< \fg_x^*
\end{CD}
$$
Since $\pi^*$ is a vector space isomorphism such that $\pi^*(\ell^*)
\subset \hat{\ell}^*$, $\hat{\phi}_V(V)$ and $\phi_V(V/\Gamma)$ are
isomorphic as subsets of vector spaces.  In particular, if $G$ is
Abelian, then $\phi_{V/\Gamma}(V/\Gamma)$ is the rational convex
polyhedral cone spanned by the orbi-weights for the action of $G_x$ on
$V/\Gamma$.  However, since $\pi^*$ need not be an isomorphism of
lattices, $\hat{\phi}_V(V) \cap \hat{\ell}^*$ and
$\phi_{V/\Gamma}(V/\Gamma) \cap \ell^*$ may not be isomorphic as
semigroups.  }
\end{remark}

As in the case of manifolds, the {\em differential} slice at $x$ is isomorphic,
as a $G_x$ space, to the product
$$
	(\fg/\fg_x)^* \times V/\Gamma .
$$
Thus, by the previous section,
a neighborhood of the $G$ orbit of $x$ in $(M, \omega)$
is equivariantly diffeomorphic to a neighborhood of the zero section
in the associated orbi-bundle
$$
	Y = G\times_{G_x} \left((\fg/\fg_x)^* \times V/\Gamma \right).
$$

\begin{Lemma}
\label{locsymp}
Let $G\cdot x$ be an {\em isotropic} orbit in $(M, \omega)$.
For every $G_x$ equivariant projection $A: \fg \to \fg_M$,
there is a symplectic form on $Y$ such that
\begin{enumerate}
\item
a neighborhood of the zero section in $Y$ is equivariantly
symplectomorphic to a neighborhood of  $G\cdot x$ in $M$, and
\item
the moment map map $\Phi_Y : Y \to \fg^*$  is given by
$$ \Phi _Y ([g, \eta, [v]])=
Ad^\dagger (g) (\eta + A^* \phi _{V/\Gamma} ([v])),$$
where $(\fg/\fg_x)^*$ is identified with the annihilator of $\fg_x$ in
$\fg^*$, $A^*: \fg^*_x \to \fg^* $ is dual to $A$,
and $\phi_{V/\Gamma} : V/\Gamma \to \fg_x^*$ is the moment
map for the action of $G_x$ on $V/\Gamma$, as in remark \ref{rem_moment}.
\end{enumerate}
\end{Lemma}
In particular, if $G$ is Abelian then
$G \cdot x$ is always isotropic and $Ad^\dagger(g)$ is trivial;
the moment map takes the form
$$
\Phi _Y ([g, \eta, [v]])= \eta + A^* \phi _{V/\Gamma} ([v]).
$$

\begin{Proof}
The construction is standard in the smooth case (cf \cite{GS}); we
simply adapt it for orbifolds.  The group $G_x$ acts on $G$ by $g_x
\cdot g = g g_x^{-1}$; this lifts to a symplectic action on $T^*G$.
The corresponding diagonal action of $G_x$ on $T^* G \times V/\Gamma$
is Hamiltonian.  An equivariant projection $A: \fg \to \fg_x$ defines
a left $G$-invariant connection 1-form on the principal $G_x$ bundle
$G \to G/G_x$, and thereby identifies $Y$ with the reduced space $(T^*
G \times V/\Gamma)_0$, thus giving $Y$ a symplectic structure.  The
$G$ moment map on $T^* G \times V/\Gamma$ descends to a moment map for
$Y$, giving the formula in $(2)$.  The proof that the neighborhoods
are equivariantly symplectomorphic reduces to a form of the
equivariant relative Darboux theorem; it is identical to the proof in
the smooth case.
\end{Proof}
\vspace{-5mm}

\begin{remark}\label{uniqueness remark}{\em
The model embedding $i:G\cdot x \hookrightarrow Y$ is {\em unique } in
the following sense.  If $i': G\cdot x \hookrightarrow (N, \sigma)$ is
any equivariant isotropic embedding into a symplectic $G$ orbifold
$(N, \sigma)$ such that the symplectic slice at $i'(x)$ is the same as
the symplectic slice $V/\Gamma$ of the model $Y$ then there exist a
neighborhood of $i(G\cdot x)$ in $Y$, a neighborhood $U'$ of
$i'(G\cdot x)$ in $N$ and a symplectic equivariant diffeomorphism
$\psi :U \to U'$ such that $i' = \psi \circ i$.  The proof of the
existence of the map $\psi$ is again, essentially, a form of the
equivariant relative Darboux theorem.}
\end{remark}
The following two lemmas are consequences of lemma  \ref{locsymp} above.

\begin{lemma}\label{lem.loc_convex}
If $G$ is a Hamiltonian torus action on a symplectic
orbifold, then the image under the moment map of a neighborhood of an
orbit is the neighborhood of a point in a rational polyhedral cone.
\end{lemma}

\begin{lemma}
If $G$ is connected, then $M^G$, the set of points which are fixed by
$G$, is a symplectic suborbifold.
\end{lemma}

We are now ready to specialize to the case of symplectic
toric orbifolds.

\begin{lemma}
\label{orbit}
Let $(M,\omega,G)$ be a symplectic toric orbifold
with moment map $\Phi_M: M \to \fg^*$.  Then for any $x \in M$,
the stabilizer of $x$ is connected.
Moreover, there is a $G$  invariant neighborhood $U$ of $G \cdot x$
on which
\begin{enumerate}
\item $\Phi_M$ induces a homeomorphism form $U/G$ to $\Phi_M(U)$.
\item the image  $\Phi_M(U)$ is the neighborhood of a point
in a simple rational polyhedral cone.
\item a neighborhood of $G \cdot x$ is classified by $\Phi_M(U)$, plus a
positive integer attached to each facet.
\end{enumerate}
\end{lemma}

\begin{proof}
Let $H$ be the stabilizer of $x$; choose a projection $A: \fg \to \fh
$.  By lemma~\ref{locsymp}, there is a neighborhood of $G \cdot x$
which is equivariantly symplectomorophic to a neighborhood of the zero
section of the model space $Y= G\times _H ( (\fg/\fh)^* \times
V/\Gamma)$, where $\Gamma$ is the orbifold structure group at $x$, and
$V/\Gamma$ is the symplectic slice at $x$.  Therefore, it suffices to
prove the above claims for the model $Y$.

Since $G$ is abelian, $H$ does not act on $(\fg/\fh)^*$.  By
assumption the action of $G$ on $M$ is effective and therefore, by
corollary~\ref{cor.loc_free} generically free.  Therefore the action
of $H$ on $V/\Gamma$ is generically free as well.

Since the action of $\Gamma$ on $V$ is generically free, it follow
that the action of $\hat{H}$, the extension $H$ by $\Gamma$, on $V$ is
generically free. Hence the symplectic representation $\hat{H}
\to Sp (V, \omega _V)$ is faithful.  That is to say we may think of
$\hat{H}$ as a compact subgroup of $Sp (2h, \R)$, where $2h=
 \dim V$, hence as a subgroup of the unitary group $U(h)$,
the maximal compact subgroup of the symplectic group.

 Also by assumption $\dim G = \frac{1}{2} \dim Y$. Consequently $\dim
H = \frac{1}{2} \dim V = h$.  The group $\tilde{H}$, the connected
component of $1$ in $\hat{H}$ is a cover of the component of $1$ in
$H$, hence is a torus of dimension $h$.  Therefore, $\tilde{H}$ is a
maximal torus of $U(h)$.  Consequently we may identify $V$ with $\C
^h$ and $\tilde{H}$ with the standard torus $T^h$ of diagonal unitary
matrices.  Since $\tilde{H}$ commutes with the action of $\Gamma$ and
since the centralizer of the maximal torus in $U(h)$ is the torus
itself, $\Gamma$ must be a subgroup of $\tilde{H}$.  In particular
{\em $\Gamma$ is abelian}.

We next argue that $H$ has to be connected.  Since $H$ is abelian, all
the elements of $H$ commute with the elements of the identity
component of $H$.  By continuity they lift to the elements of $U(h)$
that commute with the elements of $\tilde{H}$, hence by the same
argument as before, have to lie in $\tilde{H}$. Therefore $\hat{H}$ is
connected and consequently  $H$ itself is connected.

The moment maps $\hat{\phi}_V : V \to \hat\fh$ and $\phi_{V/\Gamma} :
V/\Gamma \to \fh$ are orbit maps.  The moment map map $\Phi_Y : Y \to
\fg^*$ is given by $\Phi_Y ([g, \eta, [v]])= \eta + A^* \phi
_{V/\Gamma} ([v]).$ Therefore $\Phi_Y$ is also an orbit map. Hence
the original moment map $\Phi _M$ is an orbit map.

Because the pair $(\hat{H},V)$ can be identified with $(T^h,\C^h)$, the image
$\hat{\phi}_V(V) \subset \hat{\fh}^*$ is the positive orthant.  Let $x_i
= \pi(e_i) \in \ell$, where $e_i$ is a standard basis vector in
$\hat\fh = \R^h$ and $\pi : \hat{\fh} \to \fh$ is the obvious
projection (it is actually an isomorphism).  Since
$(\pi^*)^{-1}(\hat{\phi}_V(V)) = \phi_{V/\Gamma}(V/\Gamma)$,
$\phi_{V/\Gamma}(V/\Gamma)$ is the image of the positive orthant under
$(\pi^*)^{-1}$.  An easy computation shows that $\Phi_{Y}(Y) = \cap_{i=1}^h
\{ \xi \in \fg^* \mid \left<\xi, x_i \right> \geq 0 \}$.  Therefore,
$\Phi_Y(Y)$ is a simple rational convex polyhedral cone.

For $y  = [g,\eta,[v]] \in Y$, notice that
$\Phi_Y(y)$ lies in the interior of the $i$th facet exactly
if $v_i = 0$, but the other components do not.  Here we think of $v\in V$ as
being the $h$ tuple $v= (v_1, \ldots, v_h)\in \C^h$.
In this case, it is easy to check that the structure group of $v$
is $\Z/(m_i\Z)$, where $m_i$ is the length of $x_i$.
Let $m_i$ be the number assigned to the $i$th face.

Finally, we must show that $Y$ is determined, up to equivariant
symplectomorphism, by $\phi_Y(Y)$ plus the positive integers attached
to facets.  But $\fh$ is determined by the cone $\phi_Y(Y)$, since $H$
is connected, it is also determined by $\phi_Y(Y)$.  Moreover, the
projection $\pi: \R^h \to \fh$ could be read off from the data,
because $\pi(e_i) = m_i y_i$, where $m_i$ is the integer attached to
the $i$th facet and $y_i$ is the unique primitive outward normal to
the $i$th facet.  This allows us to recover the structure group
$\Gamma$ and thereby the symplectic slice.

Since symplectically and equivariantly a neighborhood of an orbit
$G\cdot x$ in $(M, \omega, G)$ is uniquely determined by the
representation of the isotropy group of $x$ on the symplectic slice at
$x$ (cf. remark~\ref{uniqueness remark}), we are done.
\end{proof}

\section{Morse Theory}\label{section.Morse}

In this section, we extend Morse theory to orbifolds.
Since orbifolds are stratified spaces, this is simply a special case
of Morse theory on stratified spaces \cite{MTSS}.  Nevertheless, it is
an important special case which does not seem to be
readily available in literature.
We need Morse theory for the following result.

\begin{Lemma}
	\label{onedim}
$M$ be a connected compact $n$ dimensional orbifold,
and $f: M \to \R$ be a Bott-Morse function
with no critical suborbifold of index $1$ or $n-1$.
Then $M_{(a,b)} = f^{-1}(a,b)$ is connected for all $a, b \in \R$.
\end{Lemma}

We will use this result in the next section to prove that the fibers of
a torus moment map are connected, and that the image of
a compact symplectic orbifold under a torus moment map is a convex
polytope.

There are two main theorems in Morse theory; both relate the Morse
polynomial to the Poincar\'{e} polynomial.  The first states that $\cm_i
\geq \cp_i$, where $\cm_i$ is the number of critical points of index $i$,
and $\cp_i$ is the $i^{\text{th}}$ Betti number.  The second, stronger theorem
states that $\cm ( x)- \cp(x) = (1+x) Q(x)$, where $Q$ is a polynomial with
nonnegative coefficients, $\cm (x) = \sum \cm _i x^i$ is the Morse polynomial,
and $\cp(x)=\sum \cp _i x^i$ is the Poincar\'{e} polynomial.

Although the first theorem holds for orbifolds,
the second statement no longer holds;
see the following counterexample.  Luckily, most
applications of Morse theory to symplectic geometry rely only
on the first theorem.

\begin{example}
{\em
Let $M$ be torus, stood on end (see Figure~1 below).  Let $f: M \to
\R$ be the height function.  Let $\Gamma = \Z/(2\Z)$ act on $M$ by
rotating it $180$ degrees.  Then $H_0(M/\Gamma) = H_2(M/\Gamma) = \R$,
but $H_1(M/\Gamma) = 0$.  On the other hand, $\cm(x) = 1 + 2x + x^2$.}
\end{example}
%   put picture back here %%%%%%%%%%%
\setlength{\unitlength}{0.008in}
\begin{picture}(441,400)(10,-50)
\put(81,185){\ellipse{80}{156}}
\put(80,184){\ellipse{160}{240}}
\path(440,44)(439,329)
\path(441.028,321.007)(439.000,329.000)(437.028,320.993)
\path(182,183)(216,183)
\path(208.000,181.000)(216.000,183.000)(208.000,185.000)
\path(342,185)(386,185)
\path(378.000,183.000)(386.000,185.000)(378.000,187.000)
\path(322,63)	(317.538,63.540)
	(313.694,64.045)
	(310.404,64.527)
	(307.603,64.998)
	(303.210,65.951)
	(300.000,67.000)

\path(300,67)	(295.689,69.213)
	(293.183,70.748)
	(290.604,72.447)
	(285.716,75.957)
	(282.000,79.000)

\path(282,79)	(279.609,81.365)
	(276.868,84.374)
	(273.917,87.832)
	(270.899,91.542)
	(267.956,95.310)
	(265.230,98.939)
	(262.864,102.235)
	(261.000,105.000)

\path(261,105)	(258.511,109.277)
	(257.081,111.960)
	(255.625,114.812)
	(254.215,117.688)
	(252.926,120.442)
	(251.000,125.000)

\path(251,125)	(249.666,129.365)
	(248.963,132.083)
	(248.267,134.956)
	(247.599,137.835)
	(246.984,140.569)
	(246.000,145.000)

\path(246,145)	(245.010,149.381)
	(244.408,152.088)
	(243.791,154.947)
	(243.200,157.812)
	(242.677,160.539)
	(242.000,165.000)

\path(242,165)	(241.710,168.978)
	(241.516,173.834)
	(241.453,176.500)
	(241.412,179.275)
	(241.391,182.124)
	(241.390,185.009)
	(241.409,187.893)
	(241.445,190.741)
	(241.499,193.515)
	(241.569,196.179)
	(241.757,201.030)
	(242.000,205.000)

\path(242,205)	(242.481,209.221)
	(242.850,211.812)
	(243.273,214.539)
	(243.724,217.260)
	(244.179,219.837)
	(245.000,224.000)

\path(245,224)	(246.144,228.680)
	(246.906,231.551)
	(247.741,234.570)
	(248.608,237.579)
	(249.468,240.425)
	(251.000,245.000)

\path(251,245)	(252.753,249.297)
	(253.899,251.913)
	(255.148,254.653)
	(256.441,257.373)
	(257.719,259.933)
	(260.000,264.000)

\path(260,264)	(263.219,268.468)
	(265.314,271.116)
	(267.578,273.857)
	(269.894,276.553)
	(272.147,279.067)
	(276.000,283.000)

\path(276,283)	(280.956,287.218)
	(284.088,289.681)
	(287.446,292.196)
	(290.873,294.624)
	(294.213,296.826)
	(297.307,298.664)
	(300.000,300.000)

\path(300,300)	(303.061,301.052)
	(307.253,302.007)
	(309.927,302.477)
	(313.069,302.958)
	(316.739,303.462)
	(321.000,304.000)

\path(321,263)	(316.651,261.327)
	(312.917,259.833)
	(309.738,258.489)
	(307.051,257.267)
	(302.909,255.067)
	(300.000,253.000)

\path(300,253)	(296.152,248.198)
	(294.346,245.330)
	(293.000,243.000)

\path(293,243)	(290.975,238.912)
	(289.825,236.355)
	(288.660,233.641)
	(287.538,230.910)
	(286.515,228.300)
	(285.000,224.000)

\path(285,224)	(283.891,219.626)
	(283.312,216.910)
	(282.745,214.041)
	(282.210,211.166)
	(281.727,208.435)
	(281.000,204.000)

\path(281,204)	(280.412,199.831)
	(280.081,197.254)
	(279.758,194.535)
	(279.464,191.811)
	(279.225,189.221)
	(279.000,185.000)

\path(279,185)	(279.215,180.557)
	(279.454,177.832)
	(279.749,174.964)
	(280.075,172.101)
	(280.408,169.389)
	(281.000,165.000)

\path(281,165)	(281.688,160.344)
	(282.138,157.477)
	(282.644,154.458)
	(283.194,151.445)
	(283.777,148.593)
	(285.000,144.000)

\path(285,144)	(286.902,139.622)
	(288.217,137.012)
	(289.662,134.296)
	(291.148,131.606)
	(292.591,129.073)
	(295.000,125.000)

\path(295,125)	(296.375,122.700)
	(298.158,119.822)
	(302.000,115.000)

\path(302,115)	(304.652,113.231)
	(308.406,111.433)
	(310.835,110.465)
	(313.707,109.419)
	(317.077,108.272)
	(321.000,107.000)

\path(321,303)	(319.271,298.589)
	(318.090,295.306)
	(317.000,291.000)

\path(317,291)	(316.893,287.638)
	(317.140,283.449)
	(317.567,279.285)
	(318.000,276.000)

\path(318,276)	(318.908,271.452)
	(319.771,267.866)
	(321.000,263.000)

\path(321,106)	(319.271,101.589)
	(318.090,98.306)
	(317.000,94.000)

\path(317,94)	(316.893,90.638)
	(317.140,86.449)
	(317.567,82.285)
	(318.000,79.000)

\path(318,79)	(318.908,74.451)
	(319.771,70.866)
	(321.000,66.000)

\put(193,196){\makebox(0,0)[lb]{$f$}}
\put(357,196){\makebox(0,0)[lb]{$f$}}
\put(51,30){\makebox(0,0)[lb]{$\text{torus}$}}
\put(279,36){\makebox(0,0)[lb]{$\text{torus}/\Z_2$}}
\put(450,312){\makebox(0,0)[lb]{$\R$}}
\put(200,0){\makebox(0,0)[lb]{{\rm Figure 1}}}
\end{picture}

The basic definitions for  Morse theory on orbifolds are identical to their
smooth counterparts.  Let $M$ be a compact orbifold.  Let $f: M \to
\R$ be a smooth function.  We say  a critical point $x$ of
$f$ is {\em non-degenerate} if the Hessian $H(f)_x$ of $f$ is
non-degenerate.  In this case, the {\em index} of $f$ at $x$ is the
dimension of the negative eigenspace of the Hessian $H(f)_x$.
We say $f$ is {\em Morse} if  every
critical point $x$ is non-degenerate.\\[4pt]
 {\bf Notation } Denote
$f^{-1}(-\infty,a)$ by $M^-_a$ for all $a \in \R$.
\begin{Lemma}
Choose $a < b \in \R $ such that  $[a,b]$ contains no critical values.
Then $M^-_a$ is diffeomorphic to $M^-_b$,
and $f^{-1}(a)$ is diffeomorphic to $f^{-1}(b)$.
\end{Lemma}
\begin{proof}
The usual proof still applies, i.e., given a Riemannian metric on $M$,
one simply flows along the (renormalized) gradient of $f$.
\end{proof}

For critical points, the situation is only slightly more complicated.

\begin{Lemma} {\em ({\bf Morse Lemma})}
Let $p$ be a non-degenerate critical point for $f: M \to \R$ on an $n$
dimensional orbifold $M$.  There exists a neighborhood $U$ of $p$ in
$M$ and an uniformizing chart $(\tU,\Gamma,\phi)$ for $U$ such that
$\tU \subset \R^n$, $\Gamma$ acts linearly, $\phi^{-1} (p)$ is a
single point $\tilde{p}$, and
the Hessian of $f\circ \phi $ at $\tilde{p}$ equals $f\circ\phi$:
$$
H(f\circ \phi)_{\tilde{p}}  = f \circ \phi .
$$
%where $\lambda$ is the index of $f$ at $p$.
Note that the action of $\Gamma$ on $\tU$ preserves
%$\R^{\lambda}$ and $\R^{n-\lambda}$,
the positive and negative eigenspaces of the Hessian of $f\circ\phi$.
\end{Lemma}

\begin{proof}
This is simply an equivariant version of the Morse lemma for manifolds.
\end{proof}

Given a critical point $p$ of a function $f$, we can choose $\epsilon>
0$ such that $a=f(p)$ is the only critical value in
$[f(p)-\epsilon,f(p) + \epsilon]$.  Suppose further that $p$ is the
only critical point of $f$ in the level set $f^{-1}(a)$ and that the
index of $p$ is $\lambda$. Then the manifold $M^-_{a+\epsilon}$ has the
homotopy type of the space obtained by attaching the ``cell''
$D^\lambda/\Gamma$ to $M_{a -\epsilon}$ by a map from
$S^{\lambda-1}/\Gamma$ to $f^{-1}(a-\epsilon)$.  Here, $D^\lambda$ is
is the standard closed disk in the negative eigenspace of $H(f)_p$
with the action of $\Gamma$ given as above; $S^{\lambda-1}$ is its
boundary.  As before, it suffices to
check that the usual proof is equivariant.  Therefore,
$H_*(M^-_{a+\epsilon},M^-_{a-\epsilon}) = H_*(D^\lambda/\Gamma,
S^{\lambda-1}/\Gamma)$.  When we consider coefficients in $\Z$, this
can be complicated. Over $\R$, however, the following lemma is
immediate:
\begin{lemma} Let $f$ be a Morse function on an orbifold $M$. Suppose
that $p$ is the only critical point of $f$ in $f^{-1}(a-\epsilon ,a+
\epsilon )$ and that it has index $\lambda$. Then
$H_i(M^-_{a+\epsilon},M^-_{a-\epsilon};\R) = 0$ for $i \neq \lambda$;
whereas $H_\lambda(M^-_{a+\epsilon},M^-_{a-\epsilon};\R) = \R$ if
$\Gamma$ preserves the orientation of $D^\lambda$, and is trivial
otherwise.
\end{lemma}

\begin{corollary}
The number of critical points with index $i$ is greater than or equal
to the dimension of $H_i(M)$.
\end{corollary}

\begin{remark} {\em We will see in the next section that every moment map
corresponding to a circle action with isolated fixed points is a Morse
function with even indices.  Moreover, it is clear that the orbifold
structure group $\Gamma$ preserves the symplectic form, and hence the
orientation, on $D^\lambda$.  Therefore, these moment maps are perfect
Morse functions, i.e., the $i^{\text{th}}$ coefficient of the Morse
polynomial equals the $i^{\text{th}}$ coefficient of the Poincar\'{e}
polynomial: $\cm_i = \dim(H_i(M))$.}
\end{remark}

We must also consider moment maps which correspond to
circle actions with non-isolated fixed points, that is,  consider
Bott-Morse theory.
\begin{definition}
A smooth function $f: M \to \R$ is Bott-Morse if the set of critical
points is the disjoint union of suborbifolds, and if for every point
$x$ of such a suborbifold $F \subset M$, the null space of the Hessian
$H(f)_x$ is precisely the tangent space to $F$.
\end{definition}

If $f:M \to \R$ is Bott-Morse and $F$ is a critical orbifold of $f$,
the normal orbi-bundle of $F$ splits as a direct sum of vector
orbi-bundles $E^-$ and $E^+$ corresponding to the negative and
positive spectrum of the Hessian of $f$ along $F$.\footnote{Of course
one has to be careful when talking about direct sums of vector
orbi-bundles.}  Given any metric on $M$, let $D=D_F$ denote a disc
bundle of $E^-$ and $S= S_F$ denote the corresponding sphere bundle.
The index of $f$ at $F$ is the dimension of $D_F$.  In this case, we
have the following result:
\begin{lemma}
\label{morsebott}
For small $\epsilon$, $H_*(M^-_{f(F)+\epsilon},M^-_{f(F)-\epsilon}) =
H_*(D_F,S_F)$.  Moreover, the boundary map from
$H_q(M^-_{f(F)+\epsilon},M^-_{f(F)-\epsilon})$ to $H_{q-1}(M^-_{f(F) -
\epsilon})$ in the long exact sequence of relative homology is the
composition of the boundary map from $H_q(D_F,S_F)$ to $H_{q-1}(S_F)$
and the map on homology induced by the ``inclusion'' map from $S_F$ to
$M^-_{f(F) - \epsilon}$.
\end{lemma}
\begin{proof}Again, the manifold proof (see for example \cite{Chang}) can be
adapted to the case of orbifolds.
\end{proof}

We now prove lemma~\ref{onedim}.  We do it in a sequence of lemmas.

\begin{Lemma}\label{lemma2.1}
Let $F$ be an orbifold, $\pi: E \to F$ be a $\lambda$ dimensional real
vector orbi-bundle and $D(E)$ and $S(E)$ the corresponding disk and M
sphere orbi-bundles with respect to some metric. If $\lambda > 1$,
then $H_1(D(E),S(E)) = 0$.
\end{Lemma}

\begin{proof}
By the long exact sequence in relative homology it suffices to show
that the maps $H_0 (S(E)) \to H_0 (D(E))$ and $H_1 (S(E)) \to H_1
(D(E))$ induced by  inclusion are surjective.  But this follows
from two facts: the fibers of $\pi: E \to F$ are path connected and
any path in the base $F$ can be lifted to a path in the sphere bundle
$S(E)$.
\end{proof}

\begin{Lemma}
\label{Morsesurj}
Let $M$ be a connected compact orbifold, and $f: M \to \R$ be a
Bott-Morse function with no critical suborbifold of index $1$.  Then
\begin{enumerate}
\item $M^-_a = \{m \in M: f(m) <a\}$ is connected for all $a \in \R$, and
\item  if $M^-_a \neq \emptyset$, then $H_1(M^-_a) \to H_1(M)$ is
a surjection.
\end{enumerate}
\end{Lemma}

\begin{Proof}
Let $F \subset M$ be a critical suborbifold of $f$ index $\lambda$.
Let $D_F$ and $S_F$ be the disk and sphere bundles of the negative
orbi-bundle of $f$ along $F$.  Let $a = f(F)$, and let $\epsilon > 0$
be small.  We assume, for simplicity, that no other critical
suborbifold maps to $a$.

If $\lambda = 0$ then $S_F$ is empty.
Otherwise, since $\lambda \neq 1$,  $\lambda$ is greater than one
and  $H_1(D_F,S_F)$ is trivial by lemma~\ref{lemma2.1}.
In either case the map $H_1(D_F,S_F) \to H_0(S_F)$, and hence
also the map $ H_1(D_F,S_F) \to H_0(M^-_{a-\epsilon})$, is trivial.
Therefore, the following sequence is exact:
$$
0 \to H_0(M^-_{a-\epsilon}) \to H_0(M^-_{a+\epsilon}) \to
H_0(D_F,S_F) \to 0
$$
Notice that
$\dim(H_0(M^-_{a + \epsilon})) \geq \dim(H_0(M^-_{a - \epsilon}))$.  Since
$M$ is connected, this completes part (1).

If $\lambda = 0$, then $\dim(H_0(M^-_{a + \epsilon})) >
\dim(H_0(M^-_{a - \epsilon}))$.  Therefore, since $M$ is connected,
the minimum is the unique critical value of index $0$.  For any other
critical value $a$, $H_1(M^-_{a+\epsilon},M^-_{a-\epsilon}) = 0$, and
the map $H_1(M^-_{a-\epsilon}) \to H_1(M^-_{a+\epsilon})$ is a
surjection.
\end{Proof}

\begin{proof}{\bf of lemma~\ref{onedim}}
We may assume that $a$ and $b$ are regular and that $M_{(a,b)}
= \{ m \in M \mid a < f(m) < b \} \}$ is not empty.
By lemma \ref{Morsesurj},
$H_1(M^-_a) \oplus H_1(M^+_b) \to H_1(M)$ is a surjection,
where $M^+_a =  f^{-1}(a,\infty)$.
Therefore, by Mayer-Vietoris, the following sequence is exact:
$$ 0 \to H_0(M_{(a,b)}) \to H_0(M^-_a) \oplus H_0(M^+_b) \to H_0(M)
\to 0.$$ Finally, by lemma~\ref{Morsesurj}, $M^-_a$ and $M^+_b$ are
connected.
\end{proof}

\begin{remark}{\em Since $f:M \to \R$ is proper,
it follows from lemma~\ref{onedim} and a simple point set topology
argument that for any $a\in \R$ the fiber $f^{-1}(a)$ is connected.  }
\end{remark}

\section{Connectedness and Convexity}

Let $(M,\omega,T)$ be a symplectic toric orbifold with moment map
$\phi : M \to \ft^*$.  In this section we show that $\phi(M) \subset
\ft^*$ is a convex rational simple polytope.  We will prove this as a
corollary of the Atiyah-Guillemin-Sternberg convexity theorem \cite{A}
\cite{GS} for orbifolds; our proof is similar to Atiyah's. We also
prove that the fibers of the moment map $\phi$ are connected.  In
fact, as was observed by Atiyah (op.\ cit.) in the manifold case,
convexity is an easy consequence of connectedness.  In turn, the
the fibers of toral moment maps are connected because
the components of these moment maps are Bott-Morse functions with even
indices.
%The fact that the components are Bott-Morse is an easy
%consequence of the local normal form theorem.
We now give precise
statements of the main results of the section.

\begin{Theorem}
\label{connected}
Let $M$ be a Hamiltonian $T$ orbifold, $T$ a torus, with a moment map
 $\phi:M \to \ft^*$. The fibers of $\phi$ are connected.
\end{Theorem}

\begin{Theorem}
\label{convex}
Let $M$, $T$ and $\phi:M \to \ft^*$ be as in theorem~\ref{connected}
above.  Then $\Delta = \phi(M) \subset \ft^*$ is a rational convex
polytope.  In particular it is the convex hull of the image of the
points in $M$ fixed by $T$,
$$
	\phi (M) = \text{convex hull } (\phi (M^T)).
$$
\end{Theorem}

\begin{corollary}\label{cor simple}
Let $(M,\omega,T)$ be a symplectic toric orbifold with moment map
$\phi: M \to \ft^*$.  Then $\Delta = \phi(M)$ is a rational simple
polytope.
\end{corollary}
We begin the proof of the above statements with the following lemma.
\begin{lemma}\label{lemma Bott-Morse}
Let $G\times (M, \omega) \to (M, \omega )$ be a Hamiltonian group
action of a compact Lie group on a symplectic orbifold with moment map
$\phi : M \to \fg^*$.  Then for any $\xi \in \fg$ the $\xi
^{\text{th}}$ component of the moment map $\phi ^\xi := \xi \circ
\phi$ is Bott-Morse and the indices of its critical orbifolds are all
even.
\end{lemma}

\begin{proof}
This is a generalization of Theorem~5.3 of \cite{GS} and of
Lemma~(2.2) of \cite{A} to the case of orbifolds.  The proof is the
same, except one has to use the orbifold version of the equivariant
Darboux theorem (cf. lemma~\ref{locsymp} which specializes to the
equivariant Darboux theorem when the orbit is a point).
\end{proof}

\begin{Proof} {\bf \!\!of Theorem \ref{connected}}
Because the moment map $\phi$ is continuous and proper, the
connectedness of fibers is implied by the  connectedness of the
preimages of balls.  We will prove this stronger statement by induction on the
dimension of $T$.

By lemma~\ref{lemma Bott-Morse}, moment maps for circle actions are
Bott-Morse functions of even index.  Then the preimages of balls are
connected by lemma~\ref{onedim}.

Suppose now that $T$ is a $k$ dimensional torus and let $B$ be a ball
in $\ft^*$.  Let as before $\ell \subset \ft$ denote the lattice of
circle subgroups.  Then for every $0\not =\xi \in \ell$ the map $\phi
^\xi \equiv \xi \circ \phi$ is a moment map for the action of the
circle $S_\xi := \{ \exp t\xi : t\in \R\}$.  Let $\R_\xi $ denote the
set of regular values of $\phi ^\xi$.  For every $a\in \R_\xi$
the reduced space $M_{a, \xi} := (\xi \circ \phi)^{-1} (a)/S_\xi$ is a
symplectic orbifold.  The $k-1$ dimensional torus $H:= T/S_\xi$ acts
on $M_{a, \xi}$ and the action is Hamiltonian.  By inductive
assumption the preimages of balls under the $H$ moment maps $\phi ^H :
M_{a, \xi} \to \fh^*$ are connected.

The affine hyperplane $\{\eta \in \ft^*: \xi (\eta ) = a\}$ is
naturally isomorphic to the dual of the Lie algebra of $\fh$, and we
can identify the restriction of $\phi$ to $(\phi ^\xi)^{-1} (a)$ with
the pull-back of $\phi ^H$ by the orbit map $\pi :(\phi ^\xi)^{-1} (a)
\to M_{a, \xi}$.  It follows that $\phi ^{-1} (B \cap \{\eta : \xi
(\eta ) = a\}) = \pi ^{-1} ( (\phi ^H)^{-1} (B))$ is connected.

Now the set
$$
	U = \bigcup _{\xi \in \ell} \bigcup _{a\in \R_\xi } B \cap \{\eta : \xi
(\eta ) = a\}
$$
is connected and dense in the ball $B$, and its preimage  $\phi^{-1}(U)$
is connected.  Therefore the closure  $\overline{\phi^{-1} (U)}$ in $M$ is
connected.
Since $\phi $ is proper,$\overline{\phi^{-1} (U)}$ is the preimage of the
closure of $U$ in $\ft^*$, which is
$\phi ^{-1}(B)$. Hence $\phi ^{-1}(B)$ is connected.
\end{Proof}

\begin{Proof} {\bf \!\!of Theorem~\ref{convex} }
It is no loss of generality to assume that the action of $T$ is
effective, hence by corollary~\ref{cor.loc_free} free on a dense
subset.   Consequently the interior of the image $\phi (M)$ is nonempty.

To prove that $\phi (M) $ is convex it suffices to show that for any
affine line ${\cal L}\subset \ft^*$, the intersection ${\cal L}\cap
\phi (M)$ is connected.  In fact it is enough to prove this for just
rational lines, i.e. the lines of the form $\R \upsilon +a$ where
$a\in \ft^* $ and $\upsilon \in \ell ^*$, the weight lattice of $T$.

If $\upsilon \in \ell^*$ then $\ker \upsilon \subset \ft$ is the Lie
algebra of a subtorus $H= H_\upsilon$ of $T$.  A moment map $\phi ^H$
for the action of $H$ on $M$ is given by the formula $\phi ^H =i^*
\circ \phi$ where $i^*$ is the dual of the inclusion $i:\fh \to \ft$.
The fibers $(\phi ^H)^{-1} (\alpha )$ are connected by
theorem~\ref{connected}.  On the other hand
$$
 (\phi ^H)^{-1} (a) = \phi ^{-1} ((i^*)^{-1}(a)) = \phi ^{-1} (\phi
(M) \cap (a + \R \upsilon))
$$
for some $a\in (i^*)^{-1}(\alpha )$.  This proves that the image $\phi
(M)$ is convex.

Since $\phi (M)$ is compact it is, by Minkowski's theorem, the convex
hull of its extreme points.  Recall that a point $\alpha $ in the
convex set $A$ is extreme for $A$ if it {\em cannot} be written in the
form $\alpha = \lambda \beta +(1-\lambda)\gamma$ for any $\beta,
\gamma \in A$ and $\lambda \in (0,1)$. Lemma~\ref{locsymp} shows that
for any point $x$ in a Hamiltonian $T$ orbifold $M$ the image $\phi
(M)$ contains an open ball in the affine plane $\phi (x) +
\ft_x^\circ$, where $\ft_x^\circ$ is the annihilator of the isotropy
Lie algebra of $x$ in $\ft^*$.  Therefore the preimage of extreme
points of $M$ consists entirely of fixed points.

Since the set of fixed points $M^T$ is closed and $M$ is compact, $H_0
(M^T)$ is finite.  Since $\phi$ is locally constant on $M^T$, $\phi
(M^T)$ is finite.  Therefore $\phi (M) =\text{convex hull }(\phi
(M^T))$ is a convex polytope.

To prove that $\phi (M)$ is rational we need to show that the
hyperplanes supporting its facets are defined by elements of the
circle subgroup lattice $\ell \subset \ft$ (cf. the introduction).
Suppose $\xi \in \ft$ defines a hyperplane supporting a facet $F$ of
$\phi (M)$. Then the function $\phi ^\xi$ takes a global minimum on
the preimage of the facet $\phi ^{-1}(F)$.  Therefore the points in
$\phi ^{-1}(F)$ are fixed by the closure $H$ of $\{\exp t\xi: t\in \R\}$
in $T$.  It follows from the local normal form (lemma~\ref{locsymp})
that the images of the connected components of $M^H$ lie in the affine
translates of the annihilator of $\fh$ in $\ft^*$.  Since the
affine hull of the facet $F$ has codimension 1, $\fh$ is
one-dimensional, i.e. $H$ is a circle.  Therefore $\xi $ is in $\ell$.
\end{Proof}

\begin{proof}{\bf of corollary~\ref{cor simple} }
We saw in the proof of theorem~\ref{convex} that facets of the image
$\phi (M)$ correspond to one dimensional isotropy subgroups. Therefore
in order to prove that $\phi (M)$ is simple it is enough to show that
in a neighborhood of the preimage of a vertex the number of circle
isotropy groups that can occur is the same as the dimension of the
torus. But this is precisely the content of the second assertion of
lemma~\ref{orbit}.
\end{proof}

\section{From Local to Global}\label{local to global}

Let $(M,\omega,T)$ and $(M',\omega',T')$ be symplectic toric
orbifolds  with isomorphic  associated weighted polytopes.
In this section, we show that the orbifolds are isomorphic, that is,
equivariantly symplectomorphic.
The first step is to show that $M$ and $M'$ are {\em locally} isomorphic.
Then, by extending proposition~2.4 in \cite{HS} to
the symplectic category, we
show that we can ``glue'' these local isomorphisms
together to construct a global isomorphism.
Note that we may assume that $T = T'$, and that the associated
weighted polytopes are equal.

Let $T$ be a torus, and
let $(M,\omega)$ be a Hamiltonian $T$ orbifold of any dimension.
Let $\pi : M \to M/T $ be the orbit map.
Since the moment map $\phi: M \to \ft $ is $T$ invariant, we descends
to a map $\underline{\phi}:M/T \to \ft$.
Denote $\im(\phi) \subset \ft^*$ by $\Delta$.

Suppose $M'$ is another  Hamiltonian $T$ orbifold.
We say $M'$ is a orbifold {\em over $M/T$ }
if there is a continuous map  $\pi' :M' \to M/T$ which induces a
homeomorphism from $M'/T$ to $M/T$.
In this case, define
$\phi' : M' \to \ft$ by $\phi' = \underline{\phi} \circ \pi'$.

For the purposes of this section two Hamiltonian $T$ orbifolds
$(M,\omega,\pi)$ and $(M', \omega'\pi')$ over $M/T$ are {\em
isomorphic} if there exists a $T$ equivariant symplectomorphism $f: M
\to M'$ such that $\pi \circ f = \pi'$.

Two Hamiltonian $T$  orbifolds
$(M,\omega,\pi)$ and $(M', \omega',\pi')$  over $M/T$ are
{\em locally isomorphic over $\Delta$} if every point in $\Delta$ has an
open neighborhood $U$ and a
$T$ equivariant symplectomorphism $f: {\phi'}^{-1}(U) \to \phi^{-1} (U)$
such that  $\pi \circ f = \pi'$.
In this case, $\phi' :M' \to \ft$ is a moment map for the action
of $T$ on $(M',\omega')$.

\begin{remark}{\em
If $\dim T = \frac{1}{2}\dim M$, then $(M,\omega,\pi)$ and
$(M',\omega',\pi')$ are isomorphic exactly if they are equivariantly
symplectomorphic.  Furthermore, $\Delta = M/T$, so there is no need to
distinguish the two.  In contrast, if $\dim T < \frac{1}{2} \dim M$,
then $(M,\omega,\pi)$ and $(M',\omega',\pi')$ may be equivariantly
symplectomorphic but not isomorphic.  Furthermore, $\Delta \neq M/T$,
so it is important to notice that although we consider isomorphisms of
neighborhoods of fibers of the moment map, we demand that these
isomorphisms fix the orbits of $M/T$.  Since we only need the former
case, the reader may wish to assume that $\dim T = \frac{1}{2}\dim M$.
}\end{remark}

\begin{lemma}
\label{locequiv}
Let $(M,\omega,T)$ and $(M',\omega',T)$ be symplectic toric
orbifolds.  Let $\Phi: M \to \ft^*$ and $\Phi': M \to \ft^*  $
be the associated moment maps.
Assume that  $\Phi(M) = \Phi'(M')$ and that the  integers
associated to each facet are the same.
Then $M$ and $M'$ are locally isomorphic.
\end{lemma}

\begin{proof}
The proof is a direct corollary of  the local normal
form for symplectic toric orbifolds (lemma \ref{orbit}),
the connectedness of fibers of the moment map (theorem \ref{connected})
and the properness of the moment map.
\end{proof}

\begin{lemma}
\label{lem.sheaf}
Let $(M,\omega)$ be a Hamiltonian $T$ orbifold;
let $\Delta$ denote the image of its moment map.
Let $\cS$ be the sheaf over $\Delta$ defined as follows:
for each open $U \subset \Delta$,
$\cS(U)$ is the set of isomorphisms of $\phi^{-1}(U)$, that is,
the set $T$-equivariant symplectomorphisms of
$\phi^{-1}(U)$ which preserve the orbits of $T$.
Then isomorphisms classes of Hamiltonian  $T$ orbifolds
over $M/T$ which are locally isomorphic to $M$
are classified by $H^1(\Delta,\cS)$.
\end{lemma}

\begin{proof}
Let $\U = \{ U_i \}_i \in I$ be a covering of $\Delta$ such that
that there is an isomorphism $h_i : \phi^{-1}(U_i) \to {\phi'}^{-1}(U_i)$
for each $i \in I$.  Define $f_{ij}: \phi^{-1}(U_i \cap U_j) \to
\phi^{-1}(U_i \cap U_j)$ by $f_{ij} =  h_i^{-1} \circ h_j$.
These $f_{ij}$'s give a closed element of $C^1(\U,\cS)$.
Moreover, the cohomology class of this element is
independent of the choices of the isomorphisms $h_i$.

Conversely, if $\{f_{ij}\} \in C^1(\U,\cS)$ is closed, we can construct
a Hamiltonian $T$ orbifold over $M/G$ by taking the disjoint union of
the $\phi^{-1}(U_i)$'s and gluing $\phi^{-1}(U_i)$ and $\phi^{-1}(U_j)$
together using  $f_{ij}$.
\end{proof}

\begin{proposition}
\label{loctoglob}
Let $M$ and $M'$ be Hamiltonian $T$ orbifolds over $M/T$ which
are locally isomorphic.  Then $M$ and $M'$ are isomorphic.
\end{proposition}

\begin{proof}
Let $\cC^\infty$ denote the sheaf of germs of smooth functions on
$\Delta$.  Let $\underline{\ell \times \R}$ denote the sheaf of
locally constant functions with values in $\ell \times \R$.  Since
$\cC^{\infty}$ is a fine sheaf, $H^i(\Delta,C^{\infty}) = 0$ for all
$i > 0$.  Since $\Delta$ is contractable, $H^i(\Delta,\underline{\ell
\times \R}) = 0$ for all $i > 0$.

By lemma \ref{lem.sheaf} above, it suffices to show that
$H^1(\Delta,\cS) = 0.$ Therefore, by the above comments, it is
sufficient to show that that the following sequence of sheaves is
exact:
$$0 \to  \underline{\ell \times \R} \stackrel{j}{\to}
\cC^\infty \stackrel{\Lambda}{\to} \cS \to 0.$$

First we construct the map  $\Lambda: C^\infty \to \cS$.
For $U \subset \Delta$, let $f :U \to \R$ be a smooth  function.
Then the Hamiltonian vector field $X$ on $M$
of the function $f\circ \phi$ is a $T$ invariant symplectic
vector field. Moreover,
$X$ preserves $T$ orbits. Therefore $\exp(X)$ is an isomorphism
of $\phi^{-1}(U)$, i.e., $\exp(X) \in \cS(U)$.

To define $j : \underline{\R \times \ell} \to \cC^\infty$,
for any $(c,\xi) \in \R \times \ell$ and $\eta \in \ft^*$,
let $j(c,\xi)(\eta) = c +  \left<\xi,\eta\right>.$
It is clear that $j$ is injective, and that $\im(j) = \ker(\Lambda)$.

The final step is to show that $\Lambda$ is surjective.  Let $\psi$ be
an isomorphism of $\phi^{-1}(U)$, that is, a $T$-equivariant
symplectic diffeomorphism which preserves orbits.  Then obviously
$\psi$ is a $T$-equivariant diffeomorphism of $\phi^{-1}(U)$ which
preserves orbits.  Therefore, by Theorem~3.1 in \cite{HS}, there
exists a smooth $T$ invariant map $h: \phi^{-1}(U) \to T$ such that
$\psi(x) = h(x) \cdot x$.  For sufficiently small $U$, there exists a
smooth $T$ invariant map $\theta :\phi^{-1}(U) \to \ft$ such that
$\exp \circ \theta = h$.  Define a vector field $X_\theta$ on $M$ by
$X_\theta (x) = \left. \frac{d}{ds}\right |_{s= 0} \exp (s\theta
(x))\cdot x$.  A computation using a local normal form at the points
where the action of $T$ is free show that $X_\theta $ is symplectic.
Hence locally $X_\theta $ is Hamiltonian.  Choose $f$ such that $df =
i_{X} \omega$.  Since $X_theta$ is tangent to orbits of $T$, the
function $f$ Poisson commutes with all $T$ invariant functions on $M$.
Since the arguments in \cite{L} works in the case of orbifolds, it
follows that there is a smooth function $\underline{f}: \Delta \to \R$
such that $\phi^* \underline{f} = f$.  Finally, it is not hard to see
that $\Lambda(\underline{f}) = \psi$.
\end{proof}
\begin{Theorem}
Let $(M,\omega,T)$ and $(M',\omega',T')$ be symplectic toric
orbifolds with isomorphic weighted polytopes.
Then $(M,\omega, T)$ and $(M',\omega',T')$ are isomorphic.
\end{Theorem}

\begin{proof}
By lemma \ref{locequiv}, $M$ and $M'$ are locally isomorphic.
By proposition \ref{loctoglob},
locally isomorphic implies isomorphic, so we are done.
\end{proof}
\vspace{-4mm}
\begin{remark}{\em
Given any weighted polytope $\Delta$, one can construct local models
for the symplectic toric orbifold associated to $\Delta$.
Since we've shown that $H^2(\cS,\Delta) = 0$, the arguments in
\cite{HS} allow one to show that there exists a symplectic toric orbifold
which corresponds to the given weighted polytope.
However, in section \ref{surj}, we give a more explicit construction.
}\end{remark}

\section{Surjectivity}
\label{surj}

Finally,  for every weighted polytope  we construct a
corresponding  K\"{a}hler toric orbifold.
Our construction is a slight variation of Delzant's construction.

Let $\ft$ be a vector space with a lattice $\ell$.
Let $\Delta \subset \ft^*$ be a rational simplicial polytope
with $N$ facets, and a positive integer $m_i$ associated to each facet.
Then $\Delta$ can be written uniquely as
$$\Delta = \cap_{i = 1}^N
\{ \alpha \in \ft^* \mid \left<  \alpha, y_i \right> \geq \eta_i \},$$
where $y_i \in \ell$ is primitive.

Let $x_i = m_i y_i$.
Define a linear projection $\pi : \R^N \to \ft$ by $\pi(e_i) = x_i$.
Let $\fk$ be the kernel of $\pi$; let $j:\fk \to \R^N$ denote the inclusion
map.
Let $K$ be the kernel of the map from $\R^N/\Z^N $ to $\ft/\ell$
induced by $\pi$.

Let $\omega$ be the standard symplectic form on $\C^N$.
The standard action of $({S^1})^N$ on $\C^N$ has moment
map $\phi_{({S^1})^N}(z_1,\ldots,z_N) = |z_1|^2 + \cdots + |z_N|^2.$
Since $K$ is a subset of $\R^N/\Z^N$, the identification
of $\R^N/\Z^N $ with $({S^1})^N$
induces an action of $K$ on $\C^N$;
its moment map is given by $\phi_K = j^* \circ \phi_{{S^1}^N}$.

Let $(M,\sigma)$ be  the symplectic reduction of $\C^N$ by $K$ at
$j^*(m_1 \eta_1,\ldots, m_N \eta_N).$
$T = \ft/\ell = ({S^1})^N/K$ acts symplectically on $(M,\sigma)$.
Since the action of $K$ on $\C^N$ preserves
the K\"{a}hler structure on $\C^N$,  $(M,\sigma,T)$
has an induced K\"{a}hler structure.

It is easy to check that
$$
\im(\phi_T) = \cap_{i = 1}^N
\left\{ \alpha \in \ft^* \mid \left<  \alpha, m_i y_i \right> \geq m_i \eta_i
\right\}
= \Delta,
$$
where $\phi_T:M \to \ft^*$ is the moment map for $T$.

Consider $[z] \in M$, where $z = (z_1,\ldots,z_N) \in \C^N$.
It is clear that $\phi_T([z])$ lies in the interior of
the $i$th facet exactly if $\phi_{({S^1})^N}(z)$ lies in the interior of
the $i$th coordinate hyperplane, that is,  exactly if
$z_i = 0$, but $z_j \neq 0$ for $j \neq i$.
The structure group of such points is just the intersection of $K$
with the $i$th $S^1$, that is, $\Z/(m_i\Z)$.

Therefore, $(M,\sigma,T)$ is a K\"{a}hler toric variety,
and $(\Delta,\{m_i\})$ is the associated weighted polytope.

%\begin{Lemma}
%Let $\Delta$ be a rational  simple polytope.
%Let $\{m_i\}$ and $\{m'_i\}$ be two different assignments of
%integers to the facets.
%Let $(M,\sigma,T)$ and $(M',\sigma',T)$ be the associated symplectic toric
%orbifolds.
%Then there is a equivariant biholomormphic map $F : (M,T) \to (M',T)$.
%\end{Lemma}
%
%\begin{proof}
%Construct groups $K$ and $K'$ with  Lie algebras $\fK$ and $\fk'$ as above.
%Define $f: {\R^n}^* \to {\R^n}^*$ by $f(e_i) = \frac{m_i'}{m_i}$.
%Then $f$ induces a map $\hat{f}: \fk^* \to {\fk'}^*$.
%Let $\phi_K : \C^n to \fk$ and $\phi_{K'} : \C^n \to \fk'$ be
%the associated moment maps.
%\end{proof}


\begin{thebibliography}{10}
\bibitem[A]{A} M. F. Atiyah, Convexity and commuting hamiltonians,
{\it Bull. London Math. Soc.} {\bf 14} (1982), 1-15.

\bibitem[C]{Chang} K.-C. Chang,
{\em Infinite dimensional Morse theory and multiple solution problems},
Boston, Birkh\"{a}user, 1993.

\bibitem[D]{Del} T. Delzant,
Hamiltoniens periodic et images convexes de l'application moment, {\it
Bul. Soc. math. France} {\bf 116} (1988), 315-339.

\bibitem[GM]{MTSS}
 M. Goresky and R. MacPherson, {\it Stratified Morse Theory } Berlin,
New York: Springer-Verlag, 1988.

\bibitem[GS]{GS} V. Guillemin and S. Sternberg, Convexity properties of the
moment mapping I, {\it Invent. Math.} {\bf 67} (1982), 491-513.

\bibitem[HS]{HS} A. Haefliger and E. Salem, Actions of tori on orbifolds,
{\em Ann.\ Global Anal.\ Geom.} {\bf 9} (1991), 37--59.

\bibitem[K]{Kaw} K. Kawakubo, {\em The theory of Transformation Groups }
  Oxford, New York: Oxford University Press, 1991.

\bibitem[L]{L} E. Lerman, On the Centralizer of Invariant Functions
on a Hamiltonian $G$-space, {\em J.\ Diff.\ Geom.}, {\bf 30}, (1989)
805-815.

\bibitem[S]{Sa}
I. Satake, The Gauss-Bonnet theorem for V-manifolds {\em J. Math.\
Soc.\ Japan} {\bf 9} (1957), 464--492.

\end{thebibliography}
\end{document}